\documentclass[11pt,letterpaper]{article}
\usepackage[T1]{fontenc}
\usepackage{amsmath,amssymb,amsthm}
\usepackage[margin=1in]{geometry}
\usepackage{microtype}
\usepackage{xcolor}
\definecolor{paperlink}{RGB}{0,0,114}
\usepackage{hyperref}
\hypersetup{colorlinks=true,linktoc=all,linkcolor=paperlink,citecolor=paperlink,urlcolor=paperlink,pdfstartview=Fit,pdftitle={A short proof of the complete Crouzeix inequality in dimension three},pdfauthor={Sam Looi}}
\usepackage{fancyhdr}
\numberwithin{equation}{section}
\newtheorem{theorem}{Theorem}[section]
\newtheorem{lemma}[theorem]{Lemma}
\newcommand{\C}{\mathbb C}
\DeclareMathOperator{\Rea}{Re}
\newcommand{\ip}[2]{\langle #1,#2\rangle}
\title{A short proof of the complete Crouzeix inequality in dimension three}
\author{Sam Looi}
\date{31 October 2026}
\begin{document}
\maketitle

\begin{abstract}
We give a short proof of the complete Crouzeix inequality with sharp constant~$2$ for matrices of dimension at most three.
\end{abstract}

\section{Introduction}
\label{sec:introduction}

The numerical range of a matrix $A\in M_d(\C)$ is
\begin{equation*}
 W(A):=\{x^*Ax:\|x\|=1\}\,.
\end{equation*}
For a holomorphic matrix function~$F=[f_{ij}]$,  $F(A)=[f_{ij}(A)]$. We prove the following complete Crouzeix inequality for~$d\le3$, which compares the norm of a holomorphic function of a matrix with its size on the numerical range; the term ``complete'' refers to a constant independent of $m$ below.

\begin{theorem}
\label{thm:main}
Let~$1\le d\le3$ and~$A\in M_d(\C)$. For every integer~$m\ge1$ and every~$M_m(\C)$-valued function~$F$ holomorphic on a neighbourhood of~$W(A)$,
\begin{equation*}
 \|F(A)\|\le2\max_{\zeta\in W(A)}\|F(\zeta)\|\,.
\end{equation*}
\end{theorem}

The constant~$2$ is sharp already for~$d=2$ and~$m=1$. For
\begin{equation*}
 A=\begin{pmatrix}0&2\\0&0\end{pmatrix},
 \qquad f(z)=z,
\end{equation*}
we have~$\|f(A)\|=2$ and~$\max_{W(A)}|f|=1$: for~$|a|^2+|b|^2=1$, the numerical range value is~$2\overline a b$, whose modulus is at most~$1$, with equality when~$a=b=1/\sqrt2$.

Crouzeix formulated the scalar conjecture in~\cite{Crouzeix2004}, and Badea, Crouzeix and Delyon proved the complete estimate in dimension two~\cite[Section~5.2]{BCD}. A universal complete bound was obtained by Crouzeix~\cite{Crouzeix2007} and improved to~$1+\sqrt2$ by Crouzeix and Palencia~\cite{CP}. Ransford and Schwenninger gave a short proof of the abstract lemma underlying this estimate~\cite[Lemma~1.1 and Section~2]{RS}.

Scalar proofs with constant~$2$ were given by Jin~\cite[Theorem~1]{Jin2026} and Lorist and Schwenninger~\cite{LS}. \AA hag, Czy\.z and Virtanen also established the complete bound in dimension three~\cite[Theorem~1.1]{ACV}. Proofs of the complete inequality in all dimensions are now available in~\cite[Corollary~6.1]{OpenAI2026} and~\cite[Theorem~1.1]{OpenAIDirect2026}. This paper gives a short proof of the complete Crouzeix inequality in dimension three. An earlier proof of Theorem~\ref{thm:main} was completed by the author in early August 2026, and the manuscript was circulated to several colleagues throughout August. A shorter account was posted on Zenodo in early October~\cite{LooiZenodo2026}.

The proof presented here has two parts. We first prove a column estimate by completing squares and iterating over products of contractive columns. Related commuting error identities are used in~\cite[Lemma~1]{LS} and~\cite[Lemma~7.8]{ACV}. Secondly, for the complete bound, we use a Mobius variation to obtain an orthogonality relation between the components of the input and output singular vectors whenever the maximum exceeds~$2$. In dimension at most three, one of the two orthogonal spans must be one-dimensional. After transposition if needed, this reduces the problem to the column estimate. The restriction~$d\le3$ enters only in this dimension count. The column estimate is proved in Section~\ref{sec:column}, and the complete bound follows in Section~\ref{sec:completion}.

We write~$H^\infty(\Omega;M_m)$ for the bounded holomorphic~$M_m(\C)$-valued functions on~$\Omega$, with norm~$\sup_\Omega\|F\|$, and~$H^\infty(\Omega)$ in the scalar case, and we write $K\Subset\Omega$ to mean that~$K$ is compactly contained in~$\Omega$.

\section{The column estimate}
\label{sec:column}

Following the double layer construction of Delyon and Delyon \cite[Theorem~4 and Section~3]{DD} (see also~\cite[Section~2]{CP} and~\cite[Section~3]{SDV}), we use the symmetrized Cauchy formula to obtain a positive boundary measure and a commuting correction.

\begin{lemma}
\label{lem:column}
Let~$A\in M_d(\C)$, and let~$\Omega$ be a bounded convex domain with~$W(A)\Subset\Omega$. For every integer~$N\ge1$ and every~$f_1,\ldots,f_N\in H^\infty(\Omega)$ satisfying~$\sum_{j=1}^N|f_j(\zeta)|^2\le1$ for all~$\zeta\in\Omega$,
\begin{equation}
 \sum_{j=1}^N f_j(A)^*f_j(A)\le4I\,.
 \label{eq:column}
\end{equation}
\end{lemma}

\begin{proof}
We choose a positively oriented smooth convex contour~$\Gamma$ enclosing~$W(A)$ whose enclosed region is compactly contained in~$\Omega$, and write~$R_\zeta:=(\zeta I-A)^{-1}$. If~$\nu$ is the outward unit normal and~$\Rea X:=(X+X^*)/2$, the supporting line at~$\zeta$ gives
\begin{equation*}
 \Rea(\nu R_\zeta)
 =R_\zeta^*\Rea[\overline\nu(\zeta I-A)]R_\zeta\ge0\,.
\end{equation*}
Consider the positive measure
\begin{equation*}
 d\mu:=\frac1{2\pi}\Rea(\nu R_\zeta)\,ds,
 \qquad
 \ip{v}{w}_\mu:=\int_\Gamma v^*\,d\mu\,w,
 \qquad \|v\|_\mu^2:=\ip{v}{v}_\mu\,.
\end{equation*}
Since~$\nu\,ds=d\zeta/i$, we have~$\int_\Gamma d\mu=I$, so constant vectors retain their Euclidean norms. For every scalar holomorphic~$g$ on~$\Omega$, Cauchy's formula gives
\begin{equation}
 2\int_\Gamma\overline g\,d\mu=g(A)^*+E_g,
 \qquad
 E_g:=\frac1{2\pi i}\int_\Gamma\overline{g(\zeta)}R_\zeta\,d\zeta\,.
 \label{eq:completion}
\end{equation}
The resolvents, and hence~$E_g$, commute with every scalar holomorphic function of~$A$.

Suppose~$\kappa^2:=\|\sum_j C_j^*C_j\|>4$, where~$C_j:=f_j(A)$, and choose a unit vector~$x$ with
\begin{equation}
 \sum_jC_j^*C_jx=\kappa^2x\,.
 \label{eq:norming}
\end{equation}
For~$z:=\sum_j\overline{f_j}C_jx$, pointwise Cauchy--Schwarz with the square root of the density of~$\mu$ gives~$\|z\|_\mu^2\le\sum_j\|C_jx\|^2=\kappa^2$. To cancel the mixed term from~\eqref{eq:completion}, we set~$u:=z-2x$ and~$E_j:=E_{f_j}$, obtaining
\begin{equation}
 \|u\|_\mu^2+2\Rea\sum_jx^*E_jC_jx
 =\|z\|_\mu^2+4-2\kappa^2
 \le4-\kappa^2<0\,.
 \label{eq:negative}
\end{equation}

Keeping~$x$ and~$u$ fixed, for each finite column~$g=(g_\ell)$ in~$H^\infty(\Omega)$ with~$\sum_\ell|g_\ell|^2\le1$, we set~$D_\ell:=g_\ell(A)$ and
\begin{equation}
 q(g):=\|u\|_\mu^2+2\Rea\sum_\ell x^*E_{g_\ell}D_\ell x\,.
 \label{eq:sequence}
\end{equation}
We will make~$q$ arbitrarily negative by multiplying columns, contradicting a uniform lower bound. Cauchy's formula gives
\begin{equation}
 \|G(A)\|\le K_\Gamma\sup_\Gamma\|G\|,
 \qquad K_\Gamma:=\frac1{2\pi}\int_\Gamma\|R_\zeta\|\,ds
 \label{eq:cauchy-bound}
\end{equation}
for every rectangular holomorphic matrix function~$G$, independently of its size, since the integrand is~$G(\zeta)\otimes R_\zeta$. The same estimate bounds the row with entries~$E_{g_\ell}$ and the column with entries~$D_\ell$ by~$K_\Gamma$, so~$q(g)\ge\|u\|_\mu^2-2K_\Gamma^2$, independently of the column length.

The product column~$f\otimes g=(f_jg_\ell)_{j,\ell}$ is contractive because
\begin{equation}
 \sum_{j,\ell}|f_jg_\ell|^2
 =\left(\sum_j|f_j|^2\right)\left(\sum_\ell|g_\ell|^2\right)\le1\,.
 \label{eq:word-bound}
\end{equation}
By~\eqref{eq:completion} and~\eqref{eq:norming},
\begin{align}
 2\int_\Gamma\overline{g_\ell}\,d\mu\,u
 &=\sum_j\bigl((C_jD_\ell)^*+E_{f_jg_\ell}\bigr)C_jx
    -2(D_\ell^*+E_{g_\ell})x \notag\\
 &=(\kappa^2-2)D_\ell^*x
   +\left(\sum_jE_{f_jg_\ell}C_j-2E_{g_\ell}\right)x\,.
 \label{eq:prefix}
\end{align}
Multiplying by~$x^*D_\ell$, summing and using commutation, we obtain
\begin{align*}
 2q(g)-q(f\otimes g)
 &=\|u\|_\mu^2+2(\kappa^2-2)\sum_\ell\|D_\ell^*x\|^2
   -4\Rea\sum_\ell\ip{D_\ell^*x}{\overline{g_\ell}u}_\mu\\
 &\ge\sum_\ell\|\overline{g_\ell}u-2D_\ell^*x\|_\mu^2\ge0\,.
\end{align*}
The last line follows by completing squares, using~$\kappa^2>4$ and~$\sum_\ell\|\overline{g_\ell}u\|_\mu^2\le\|u\|_\mu^2$. Since~$q(f)<0$ by~\eqref{eq:negative}, iterating~$q(f\otimes g)\le2q(g)$ from~$g=f$ makes~$q$ unbounded below, a contradiction.
\end{proof}

\section{The complete bound}
\label{sec:completion}

We use a M\"obius variation from Hartz--McCarthy~\cite[Lemma~4.1]{HM} to show that, when the extremal norm exceeds~$2$, the input and output component spans are orthogonal.

\begin{proof}[Proof of Theorem~\ref{thm:main}]
We fix a bounded convex domain~$\Omega$ with~$W(A)\Subset\Omega$ and an integer~$m\ge1$. The supremum~$c:=\sup_G\|G(A)\|$ over the unit ball of~$H^\infty(\Omega;M_m)$ is finite by~\eqref{eq:cauchy-bound}. By Montel's theorem, a maximizing sequence has a subsequence converging locally uniformly to a function in the same unit ball. By Cauchy's formula on a fixed interior contour, we have convergence of the evaluations at~$A$, so the limit attains~$c$. If~$c>2$, we choose an extremal~$F$ and unit singular vectors $\xi, \eta$ of $T:= F(A)$ satisfying
\begin{equation}
 T\xi=c\eta,
 \qquad T^*\eta=c\xi\,.
 \label{eq:singular}
\end{equation}
For unitary~$U\in M_m(\C)$ and real~$t$ with~$|t|<1$, we set~$F_t:=(F+tU)(I+tU^*F)^{-1}$. The inverse exists by the Neumann series, and~$F_t$ remains contractive because
\begin{equation*}
 I-F_t^*F_t
 =(1-t^2)(I+tF^*U)^{-1}(I-F^*F)(I+tU^*F)^{-1}\ge0\,.
\end{equation*}
Thus~$\Rea(\eta^*F_t(A)\xi)\le c$, with equality at~$t=0$. Differentiating there, using~\eqref{eq:singular} and writing~$\widehat U:=U\otimes I_d$, gives
\begin{equation*}
 0=\Rea\bigl[\eta^*(\widehat U-T\widehat U^*T)\xi\bigr]
  =(1-c^2)\Rea(\eta^*\widehat U\xi)\,.
\end{equation*}
Recalling that $U$ was an arbitrary unitary, we may replace $U$ by~$iU$, which gives~$\eta^*\widehat U\xi=0$. We write~$\xi=\sum_j e_j\otimes x_j$ and~$\eta=\sum_i e_i\otimes y_i$. Every matrix unit is half the difference of two unitaries (a suitable permutation matrix and the matrix obtained by negating one of its columns), so
\begin{equation}
 y_i^*x_j=0\qquad\text{for every }i,j\,.
 \label{eq:orthogonality}
\end{equation}
The nonzero spans of the~$x_j$ and~$y_i$ are orthogonal in~$\C^d$. Since~$d\le3$, one has dimension one.

Transposition preserves~$W(A)$ and gives~$F^{\mathsf T}(A^{\mathsf T})=T^{\mathsf T}$ with singular pair~$(\overline\eta,\overline\xi)$. Transposing if necessary, we may write~$\xi=\alpha\otimes x$ with~$\|\alpha\|=\|x\|=1$. The contractive column~$f:=F\alpha=(f_j)_j$ then satisfies, by Lemma~\ref{lem:column},
\begin{equation*}
 c^2=\|T(\alpha\otimes x)\|^2
     =\sum_j\|f_j(A)x\|^2\le4\,,
\end{equation*}
a contradiction.

Homogeneity gives~$\|F(A)\|\le2\sup_\Omega\|F\|$. For~$F$ holomorphic near~$W(A)$, compactness and convexity give bounded convex domains~$\Omega_\varepsilon:=\{\zeta:\operatorname{dist}(\zeta,W(A))<\varepsilon\}$ whose closures lie in the domain of~$F$ for small~$\varepsilon>0$. As~$\varepsilon\to0^+$, continuity gives~$\sup_{\Omega_\varepsilon}\|F\|\to\max_{W(A)}\|F\|$.
\end{proof}

\medskip
\noindent California Institute of Technology, Pasadena, CA 91125, USA\\
\href{mailto:s.looi9309@gmail.com}{Email}
\end{document}